\input amstex
\documentstyle{amsppt}\nologo\footline={}
\hsize450pt

\def\SO{\mathop{\text{\rm SO}}}
\def\tr{\mathop{\text{\rm tr}}}

\topmatter
\title A couple of real hyperbolic disc bundles over surfaces\endtitle
\author Sasha Anan$'$in, Philipy V.~Chiovetto\endauthor
\thanks The second author is supported by the grant Inicia\c c\~ao Cient\'\i fica FAPESP 2014/26282-5\endthanks
\address Departamento de Matem\'atica, ICMC, Universidade de S\~ao Paulo, Caixa Postal 668,\newline13566-590--S\~ao
Carlos--SP, Brasil\endaddress
\email sasha\_a\@icmc.usp.br\endemail
\address Departamento de Matem\'atica, ICMC, Universidade de S\~ao Paulo, Caixa Postal 668,\newline13566-590--S\~ao
Carlos--SP, Brasil\endaddress
\email philipy.chiovetto\@usp.br\endemail
\abstract
Applying the techniques developed in [AGG], we construct new real hyperbolic manifolds whose underlying topology is that of
a disc bundle over a closed orientable surface. By the Gromov-Lawson-Thurston conjecture [GLT], such bundles $M\to S$
should satisfy the inequality $|eM/\chi S|\leqslant1$, where $eM$ stands for the Euler number of the bundle and $\chi S$,
for the Euler characteristic of the surface. In~this paper, we construct new examples that provide a maximal value of
$|eM/\chi S|=\frac35$ among all known examples. The former maximum, belonging to Feng Luo [Luo], was $|eM/\chi S|=\frac12$.
\endabstract
\endtopmatter
\document

{\hfill\it Dedicated to Krolik}

\bigskip

\centerline{\bf1.~Introduction}

\medskip

Topologically or differentially, every open disc bundle $M\to S$ over a closed connected orientable surface $S$ can be
completely characterized by two numbers: the Euler characteristic $\chi S$ of the surface and the Euler number $eM$ of the
bundle, i.e., the number of self-intersections of a section of the bundle. Note that, taking an unramified finite cover of
$S$ and pullbacking the bundle, one gets the same value of $|eM/\chi S|$.

The conjecture of Gromov, Lawson, and Thurston [GLT, p.~28] suggests a numerical criterion deciding whether a bundle can
be equipped with a complete real hyperbolic geometry.

\medskip

{\bf1.1.~GLT-conjecture.} A disc bundle $M\to S$ over a closed connected orientable surface $S$ of genus $g\geqslant2$
admits a complete real hyperbolic structure iff $|eM/\chi S|\leqslant1$.

\medskip

In [AGG], Conjecture 1.1 was extended (with the same bound) to the complex hyperbolic case.

\medskip

{\bf1.2.~Known results.} The best proven upper bound belongs to Misha Kapovich [Kap] who showed that
$|eM|\leqslant\exp\Big(\exp\big(10^8|\chi S|\big)\Big)$ for any complete real hyperbolic $4$-manifold homotopically
equivalent to a closed orientable surface; so, without actually using the fact that it admits a disc bundle structure. (In
this case, $eM$ stands for the self-intersection of the generator of $H_2M$ represented by a homotopy equivalence
$S\to M$.) It is worthwhile mentioning that, in such settings, Nicholaas H.~Kuiper [Kui2] constructed examples with
$|eM/\chi S|>\frac2{\sqrt3}>1$.

In the other direction, the best results belong to N.~H.~Kuiper [Kui1] and to Feng Luo [Luo]. In~[Kui1, Theorem 6, p.~68],
it was constructed a series of disc bundles admitting complete real hyperbolic\break
geometry with any rational value of $|eM/\chi S|$ in the interval $\big[0,\frac13\big]$. (Though, we are not sure that this
result is literally correct as there are a few miscalculations in the exposition.) F.~Luo constructed an example with a
maximal known (before our paper) value $|eM/\chi S|=\frac12$. Since the surface $S$ in F.~Luo's example has genus $2$,
taking a finite unramified cover of $S$, one gets examples satisfying the relation $|eM/\chi S|=\frac12$ with $S$ of an
arbitrary genus $g\geqslant2$.

\medskip

{\bf1.3.~Main result.} Using the ideas of [AGG] and a coordinate free approach to hyperbolic geometry exposed in [AGr2,
AGr3], we construct $3$ new examples of disc bundles with $|eM/\chi S|=\frac35$.

Two of them have $eM=12$ with $S$ of genus $11$. The third one has $eM=24$ with $S$ of genus $21$. The first two come from
right-angled necklace polyhedra with $84$ codimension $1$ faces and $21$ cycles of codimension $2$ faces. The third one has
a similar flavor with $164$ codimension $1$ faces and $41$ cycles of codimension $2$ faces.

In fact, we construct $3$ orbifolds whose $4$-sheeted covers provide the above manifolds. Each orbifold comes from a
right-angled necklace polyhedron $P$ symmetric with respect to a regular elliptic isometry $r$ of order $24$ for the first
$2$ examples and of order $44$ for the third one. The face-pairing of these polyhedra are reflections in their totally
geodesic codimension $1$ faces.

For other examples of disc bundles of a similar type, see Theorem 3.9.

\bigskip

\centerline{\bf2.~Construction}

\medskip

We start the construction by fixing a regular elliptic isometry $r$ of order $n$. Denote by $P_1$ and $P_2$ its totally geodesic
$r$-stable planes. The planes intersect orthogonally at the unique $r$-fixed point $b$.

Next, we pick a generic totally geodesic hyperplane $H_0\subset\overline{\Bbb H}_\Bbb R^4$ in the real hyperbolic $4$-space
$\overline{\Bbb H}_\Bbb R^4$ considered together with its absolute (see the beginning of Section 3). Geometrically, the pair
$r,H_0$ is given by the distances from $b$ to the intersections $H_0\cap P_1$ and $H_0\cap P_2$. So, it can be described by
means of $2$ convenient real parameters $x_1,x_2$ responsible for these distances.

Then we copy the hyperplanes, $H_i:=r^iH_0$. The conditions that
$C_i:=H_i\cap H_{i+1}\not\subset\partial\overline{\Bbb H}_\Bbb R^4$ and that the other pairs of hyperplanes $H_i$ and
$H_j$, $i-j\not\equiv_n\pm1$, are ultraparallel, i.e., $H_i\cap H_j=\varnothing$, is equivalent by Lemma 3.1 to a finite
number of inequalities linear in $x_1,x_2$. Thus, we arrive at a convex region $R\subset\Bbb R^2(x_1,x_2)$. In what
follows, we assume $(x_1,x_2)\in R$.

Intersecting those closed half-spaces limited by the $H_i$'s that contain the point $b$, we get a convex polyhedron $P$
bounded by closed solid cylinders $B_i\subset H_i$ and by a piece $\partial_1P\subset\partial\overline{\Bbb H}_\Bbb R^4$ of
the absolute bounded by a torus $T\subset\partial\overline{\Bbb H}_\Bbb R^4$. In turn, the solid cylinder $B_i$ is bounded
inside $H_i$ by the ultraparallel totally geodesic planes $C_{i-1},C_i$ and by a cylinder in $T$; such cylinders form the torus
$T$. Denoting $\partial_0P:=\bigcup\limits_{i=1}^nB_i$, we see that $\partial_0P$ is a solid torus bounded by $T$ and that
$\partial P=\partial_0P\sqcup_T\partial_1P$.

Every solid cylinder $B_i$ is fibred by totally geodesic planes called {\it slices\/} of $B_i$. Indeed, the geodesic
segment $\Gamma_i$ that joins the closest points in $C_{i-1}$ and in $C_i$ lists the fibres in question: through any
$p\in\Gamma_i$, we have, inside $H_i$, a totally geodesic plane orthogonal to $\Gamma_i$. Note that $C_{i-1}$ and $C_i$ are
among the slices; they are the {\it initial\/} and the {\it final\/} slices. Denote by $M_i$ the {\it middle slice,} i.e.,
the one passing through the middle point of $\Gamma_i$.

Every solid cylinder $B_i$ is also fibred by closed segments called {\it strings\/} of $B_i$. They can be described as
follows. Take any totally geodesic plane $F$ such that $\Gamma_i\subset F\subset H_i$. The intersection $F\cap B_i$ is bounded
in $F$ by the geodesics $F\cap C_{i-1}$ and $F\cap C_i$, both orthogonal to $\Gamma_i$, and by two arcs on the absolute. A {\it
string\/} of $B_i$ is the segment of a line, inside some $F$, equidistant from the geodesic containing the segment $\Gamma_i$.
The segment $\Gamma_i$ and the mentioned two arcs on the absolute are among the strings of $B_i$. Clearly, the reflection
$\sigma_i$ in the middle slice $M_i$ of $B_i$ stabilizes any string of $B_i$ and interchanges the endpoints of the string.

Pick a point $q_0\in\partial\overline{\Bbb H}_\Bbb R^4\cap C_0$ and let $q_0\in s_n\subset B_n$ be the string of $B_n$ that
contains $q_0$. Then $q_{n-1}:=\sigma_nq_0\in\partial\overline{\Bbb H}_\Bbb R^4\cap C_{n-1}$ is the other endpoint of
$s_n$. Next, we take the string $s_{n-1}$ of $B_{n-1}$ such that $q_{n-1}\in s_{n-1}\subset B_{n-1}$, and so on. Finally,
we get a simple curve $s:=s_1\cup s_2\cup\dots\cup s_n$ with the endpoints
$q_0\in\partial\overline{\Bbb H}_\Bbb R^4\cap C_0$ and
$q'_0:=\sigma_1\sigma_2\dots\sigma_nq_0\in\partial\overline{\Bbb H}_\Bbb R^4\cap C_0$, where
$\partial\overline{\Bbb H}_\Bbb R^4\cap B_i\supset s_i$ is a string of $B_i$ for all $1\leqslant i\leqslant n$. We call
such a curve $s\subset T$ the {\it string\/} of $P$ {\it generated\/} by
$q_0\in\partial\overline{\Bbb H}_\Bbb R^4\cap C_0$.

\medskip

{\bf2.1.~Lemma {\rm(cf.~[AGG, Lemma 2.25, p.~4317])}.} {\sl If\/ $H_0$ is orthogonal to\/ $P_1$ {\rm(}in terms of the
parameters, this means that\/ $x_2=0${\rm),} then any string of\/ $P$ is closed and contractible in\/ $\partial_1P$. So,
$\partial_1P$ is a solid torus and the slice bundle of\/ $\partial_0P$ is extendable to\/ $P$ in this case.}

\medskip

{\bf Proof.} When $H_0$ is orthogonal to $P_1$, the intersection $Q:=P\cap P_1$ is a regular $n$-gon centred at $b$ in the
hyperbolic plane $P_1$. The polyhedron $P$ is simply the union of all those totally geodesic planes orthogonal to $P_1$
that pass through a point of $Q$. The slices of $\partial_0P$ are built over the points of the boundary $\partial Q$. Thus,
we get the slice bundle of $\partial_0P$ extended to $P$.

The isometry $\sigma_i$ is the trivial extension of the reflection in the middle point of the corresponding side of $Q$.
Hence, the isometry $\sigma_1\dots\sigma_n$ is a trivial extension of an elliptic isometry of $P_1$ with the fixed point
$C_0\cap P_1$. In other words, the restriction $\sigma_1\dots\sigma_n|_{C_0}$ is the identity, implying that any string of
$P$ is closed.

In order to visualize a contraction of a closed string $s\subset T$ in $\partial_1P$, one can simply shrink the $n$-gon~$Q$
(say, keeping its $r$-rotational symmetry about $b$).

The fact that $\partial_1P$ is a solid torus follows from the Dehn lemma
$_\blacksquare$

\medskip

If any string of the polyhedron $P$ is closed, we say that $P$ is {\it fibred.} As we saw, this is equivalent to
$\sigma_1\dots\sigma_n|_{C_0}=1_{C_0}$. Denote $\sigma:=\sigma_0$. Then $\sigma_i=r^i\sigma r^{-i}$ and
$\sigma_1\dots\sigma_n=(r\sigma)^n$ because $r^n=1$. Since $C_0$ is $r\sigma$-stable, we conclude that $P$ is fibred iff
$r\sigma|_{C_0}$ is an elliptic isometry of $C_0$ whose order divides $n$.

It follows from the connectedness of the region $R$ and from Lemmas 2.1 and 3.1 that $\partial_1P$ is a solid torus. By
[AGG, Lemma 2.19, p.~4312], $P$ is topologically a closed $4$-ball and the slice bundle of $\partial_0P$ is extendable to
$P$.

Suppose that $P$ is fibred. Take a couple of disjoint strings $s,s'\subset T$ of $P$. Since $s$ and $s'$ are contractible
inside $P$ with $2$-discs $D,D'\subset P$, $s=\partial D$ and $s'=\partial D'$, one can calculate the algebraic number of
intersections of $D$ and $D'$. This number $eP$ (clearly independent of the choice of $s,s'$) is the {\it Euler number\/}
of the fibred polyhedron $P$.

Suppose that a fibred polyhedron $P$ is equipped with face-pairing isometries for the codimension $1$ faces $B_i$ and that the
conditions of Poincar\'e's polyhedron theorem (PPT) are satisfied (see, for example, [AGr1]). Then we obtain a disc bundle
$M$ over a $2$-dimensional orbifold $S$ and $eP$ is the Euler number of this bundle because the (face-pairing) isometries
preserve the slice bundles and the string bundles of the~$B_i$'s. (See, for instance, [BoS] for a treatment of bundles over
orbifolds. Another option is to glue a few copies of $P$ forming a fundamental polyhedron for a manifold and to note that
$eP$ and $\chi S$ get multiplied by the number of copies.)

A closed curve $c\subset T$ that generates the group $H_1(\partial_0P,\Bbb Z)$ is said to be {\it trivializing\/} if
$[c]=0$ in $H_1(\partial_1P,\Bbb Z)$. In the case considered in Lemma 2.1, any string of $P$ is trivializing.

\medskip

{\bf2.2.~Lemma {\rm[AGG, Remark 2.22, p.~4314]}.} {\sl Let\/ $P$ be a fibred polyhedron, let\/ $T\supset s$ be a string
of\/~$P$, and let\/ $T\supset c$ be a trivializing curve of\/ $P$. Then\/ $eP=\#s\cap c$. In other words, $[s]=eP\cdot[g]$
in the group\/ $H_1(\partial_1P,\Bbb Z)$, where\/ $g:=\partial\overline{\Bbb H}_\Bbb R^4\cap C_0$}
$_\blacksquare$

\medskip

Now we get a tool to measure the Euler number $eM$. Given fibred polyhedron satisfying the conditions of PPT, one can
deform it into a `plane' one (dealt with in Lemma 2.1) because the region $R$ is convex. At~the beginning of the
deformation, the Euler number $eP$ of a `plane' polyhedron equals $0$ by Lemmas~2.1 and 2.2. During the deformation, we
keep track of how many times the polyhedron becomes fibred, i.e., how many times the string of the polyhedron becomes
closed. Of course, counting these events, we~should take care of the signs. So, it is better to say that the Euler number
of the fibred polyhedron at the end of the deformation equals the algebraic number of times it was fibred during the
deformation, including the last moment and not including the initial `plane' moment.

We apply this method at the end of the proof of Theorem 3.9. The count is simple there because the chosen deformation
provides a monotonic evolution of a string.

\bigskip

\centerline{\bf3.~Calculation}

\medskip

Let $b_1,b_2,b_3,b_4,b$ be an orthonormal basis of signature $----+$ in an $\Bbb R$-linear space $V$ equipped with a
symmetric bilinear form $\langle-,-\rangle$. Then the real hyperbolic space $\Bbb H_\Bbb R^4$, its absolute
$\partial\Bbb H_\Bbb R^4$, and $\overline{\Bbb H}_\Bbb R^4:={\Bbb H}_\Bbb R^4\sqcup\partial{\Bbb H}_\Bbb R^4$ are known to
be identified respectively with
$$\Bbb H_\Bbb R^4:=\big\{p\in\Bbb P_\Bbb RV\mid\langle p,p\rangle>0\big\},\qquad\partial\overline{\Bbb H}_\Bbb
R^4:=\big\{p\in\Bbb P_\Bbb RV\mid\langle p,p\rangle=0\big\},\qquad\overline{\Bbb H}_\Bbb R^4:=\big\{p\in\Bbb
P_\Bbb RV\mid\langle p,p\rangle\geqslant0\big\}.$$

Pick some numbers $k,m,n\in\Bbb N$ such that $n$ is even and $1<k<m<\frac n2$. Denote by
$$r:=\left[\smallmatrix
c_1&-s_1&0&0&0\\s_1&c_1&0&0&0\\0&0&c_m&-s_m&0\\0&0&s_m&c_m&0\\0&0&0&0&1\endsmallmatrix\right],\qquad
p_0:=\left[\smallmatrix\sqrt{x_1}\\0\\\sqrt{x_2}\\0\\\sqrt{x_1+x_2-1}\endsmallmatrix\right],\qquad0<x_1,\quad0\leqslant
x_2,\quad1<x_1+x_2,$$
the regular elliptic isometry of $\Bbb H_\Bbb R^4$ with the unique fixed point $b\in\Bbb H_\Bbb R^4$ and a point
$p_0\in\Bbb P_\Bbb RV\setminus\overline{\Bbb H}_\Bbb R^4$, both written in the above basis, where $x_1,x_2$ are some real
parameters subject to the displayed inequalities, $c_i:=\cos\frac{2i\pi}n$, and $s_i:=\sin\frac{2i\pi}n$. Clearly, $r^n=1$,
$r\in\SO V$, and $\langle p_0,p_0\rangle=-1$. For any $i\in\Bbb Z$, denote
$$p_i:=r^ip_0=\left[\smallmatrix
c_i&-s_i&0&0&0\\s_i&c_i&0&0&0\\0&0&c_{mi}&-s_{mi}&0\\0&0&s_{mi}&c_{mi}&0\\0&0&0&0&1\endsmallmatrix\right]\cdot
\left[\smallmatrix\sqrt{x_1}\\0\\\sqrt{x_2}\\0\\\sqrt{x_1+x_2-1}\endsmallmatrix\right]=\left[\smallmatrix
c_i\sqrt{x_1}\\s_i\sqrt{x_1}\\c_{mi}\sqrt{x_2}\\s_{mi}\sqrt{x_2}\\\sqrt{x_1+x_2-1}\endsmallmatrix\right].$$
Obviously,
$$g_i:=\langle p_0,p_i\rangle=(1-c_i)x_1+(1-c_{mi})x_2-1.$$
Denote by $H_i\subset\overline{\Bbb H}_\Bbb R^4$ the totally geodesic hyperplane corresponding to the $\Bbb R$-linear
subspace $p_i^\perp\leqslant V$. Let $C_i:=H_i\cap H_{i+1}$.

\medskip

{\bf3.1.~Lemma.} {\sl The conditions that\/ $C_i\not\subset\partial\overline{\Bbb H}_\Bbb R^4$ and that\/
$H_i\cap H_j=\varnothing$ for all\/ $i,j$ such that\/ $i-j\not\equiv_n\pm1$ are equivalent to the inequalities
$$0<x_1,\qquad0\leqslant x_2,\qquad1<x_1+x_2,\qquad(1-c_1)x_1+(1-c_m)x_2<2<(1-c_i)x_1+(1-c_{mi})x_2,\leqno{\bold{(3.2)}}$$
for\/ $2\leqslant i\leqslant\frac n2$. The convex region\/ $R\subset\Bbb R^2(x_1,x_2)$ given by the inequalities\/
{\rm(3.2)} is nonempty because it contains the segment\/ $A:=\big(\frac2{1-c_2},\frac2{1-c_1}\big)\times0\subset R$.}

\medskip

{\bf Proof.} The condition $C_i\not\subset\partial\overline{\Bbb H}_\Bbb R^4$ means that the form on $\Bbb Rp_i+\Bbb Rp_{i+1}$
is of signature $--$, i.e., that
$\det\left[\smallmatrix\langle p_i,p_i\rangle&\langle p_i,p_{i+1}\rangle\\\langle p_{i+1},p_i\rangle&\langle
p_{i+1},p_{i+1}\rangle\endsmallmatrix\right]>0$
by Sylvester's criterion. Taking into account that $\langle p_k,p_k\rangle=-1$ for all~$k$, it is equivalent to
$-1<\langle p_i,p_{i+1}\rangle<1$. Similarly, the condition that $H_i\cap H_j=\varnothing$ is equivalent to
$\langle p_i,p_j\rangle\notin[-1,1]$ because it means that the signature of the form on $\Bbb Rp_i+\Bbb Rp_j$ is $-+$. Since
$0\leqslant g_i+1$ and $0<g_1+1$, due to the symmetry related to the action of $r$, we arrive at the inequalities (3.2). For
$x_2=0$, the inequalities take the form $1<x_1$ and $(1-c_1)x_1<2<(1-c_2)x_1$, i.e., the form
$1<\frac2{1-c_2}<x_1<\frac2{1-c_1}$ with $-1<c_2<c_1$
$_\blacksquare$

\medskip

In the sequel, we assume $(x_1,x_2)\in R$.

Obviously, $C_i\cap C_j=\varnothing$ unless $i\equiv_nj$. Therefore, the ultraparallel planes
$C_{i-1},C_i\subset H_i\simeq\overline{\Bbb H}_\Bbb R^3$ limit in $H_i$ a solid cylinder $B_i$, and these cylinders form a
solid torus $\partial_0P:=\bigcup\limits_{i=1}^nB_i$ bounded by a torus $T\subset\partial\overline{\Bbb H}_\Bbb R^4$ as
claimed in Section 2.

All the $C_i$'s are in the same closed half-space of $\overline{\Bbb H}_\Bbb R^4$ limited by a hyperplane $H_j$ as, otherwise,
there would exist some $C_{i-1}$ and $C_i$, both disjoint with $H_j$, in different half-spaces, which would cause an impossible
intersection $H_i\cap H_j\ne\varnothing$.

Let $P$ stand for the intersection of those closed half-spaces limited by the $H_i$'s that contain the point~$b$; the
$B_i$'s are codimension $1$ faces of $P$ and the $C_i$'s are codimension $2$ faces of $P$. By Lemma~2.1, $\partial_1P$ is a
solid torus and $\partial P=\partial_0P\sqcup_T\partial_1P$ is a closed $4$-ball.

Denote by $\sigma$ the reflexion in the middle slice $M_0$ of $B_0$ and by $\tau$, the reflection in $H_0$. Clearly,
$\sigma_i:=r^i\sigma r^{-i}$ is the reflection in the middle slice $M_i$ of $B_i$ and $\tau_i:=r^i\tau r^{-i}$ is the
reflection in $H_i$.

\medskip

{\bf3.3.~Lemma.} {\sl Suppose that\/ $g_1=0$, i.e.,
$$(1-c_1)x_1+(1-c_m)x_2=1.\leqno{\bold{(3.4)}}$$
Then\/ $(\tau_{i+1}\tau_i)^2=1$ for all\/ $i$ and the polyhedron\/ $\Bbb H_\Bbb R^4\cap P$ endowed with the face-pairing
isometries\/ $\tau_i$ {\rm(}identifying every codimension\/ $1$ face\/ $B_i$ with itself\/{\rm)} satisfies the conditions
of Poincar\'e's polyhedron theorem.}

\medskip

{\bf Proof.} The hyperplanes $H_i$ and $H_{i+1}$ are orthogonal along $C_i$ because $\langle p_i,p_{i+1}\rangle=g_1=0$. As
the reflection $\tau_i$ is given by the rule $\tau_i:v\mapsto v+2\langle v,p_i\rangle p_i$, the equality
$(\tau_{i+1}\tau_i)^2=1$ follows straightforwardly from $\langle p_i,p_{i+1}\rangle=0$ and
$\langle p_i,p_i\rangle=\langle p_{i+1},p_{i+1}\rangle=-1$.

Since $\tau_i$ sends the interior of $P$ into the exterior of $P$ and every geometric cycle of codimension $2$ faces of $P$
has length $4$ and, therefore, total angle $2\pi$, the conditions of PPT are satisfied (see [AGr1, Theorem 3.2, p.~303] and
[AGr1, Proposition 2.1, p.~300])
$_\blacksquare$

\medskip

{\bf3.5.~Lemma.} {\sl Let\/ $U:=\Bbb Rp_0+\Bbb Rp_1$ and\/ $W:=\Bbb Rp_0+\Bbb R(p_1-p_{n-1})$. Then\/ $C_0$ and\/ $M_0$
correspond respectively to the\/ $\Bbb R$-linear subspaces\/ $U^\perp$ and\/ $W^\perp$. The isometries\/ $\sigma$ and\/
$r\sigma$ are given by the rules
$$v\mapsto v+2\langle v,p_0\rangle p_0+\frac{\langle v,p_1-p_{n-1}\rangle}{g_2+1}(p_1-p_{n-1}),\leqno{\bold{(3.6)}}$$
$$r\sigma v=rv+2\langle v,p_0\rangle p_1+\frac{\langle v,p_1-p_{n-1}\rangle}{g_2+1}(p_2-p_0).\leqno{\bold{(3.7)}}$$
The\/ $\Bbb R$-linear subspace\/ $U$ is\/ $r\sigma$-stable and\/ $\left[\smallmatrix0&1\\-1&2g_1\endsmallmatrix\right]$ is
the matrix of\/ $r\sigma|_U$ in the basis\/ $p_0,p_1$. The point\/
$f_0:=(1-g_1)b+\langle b,p_0\rangle(p_0+p_1)\in\Bbb H_\Bbb R^4\cap C_0$ is a fixed point of\/ $r\sigma$.}

\medskip

{\bf Proof.} By definition, $C_0$ corresponds to $U^\perp$. In other words, $C_0$ corresponds to
$p_0^\perp\cap(p_1+g_1p_0)^\perp$~and, similarly, $C_{n-1}$ corresponds to $p_0^\perp\cap (p_{n-1}+g_1p_0)^\perp$, where
$p_1+g_1p_0,p_{n-1}+g_1p_0\in p_0^\perp$. Since
$$\langle p_1+g_1p_0,p_1+g_1p_0\rangle=\langle p_1+g_1p_0,p_1\rangle=-1+g_1^2<0,$$
$$\langle p_{n-1}+g_1p_0,p_{n-1}+g_1p_0\rangle=\langle p_{n-1}+g_1p_0,p_{n-1}\rangle=-1+g_1^2<0,$$
$$\langle p_1+g_1p_0,p_{n-1}+g_1p_0\rangle=\langle p_1+g_1p_0,p_{n-1}\rangle=g_2+g_1^2>0,$$
the middle point of $\Gamma_0$ equals $m_0:=(p_1+g_1p_0)+(p_{n-1}+g_1p_0)=p_1+p_{n-1}+2g_1p_0$ and $M_0$ corresponds to
$p_0^\perp\cap\big((p_1+g_1p_0)-(p_{n-1}+g_1p_0)\big)^\perp$.

Clearly, $p_0$ and $p_1-p_{n-1}$ are orthogonal. Hence, the rule (3.6) acting as $v\mapsto v$ for any $v\in W^\perp$ and as
$v\mapsto-v$ for any $v\in W$ in view of $\langle p_1-p_{n-1},p_1-p_{n-1}\rangle=-2g_2-2$ defines the isometry~$\sigma$.
The~formula~(3.7) is now immediate. It implies the equalities $r\sigma p_0=-p_1$ and $r\sigma p_1=p_0+2g_1p_1$ providing
the indicated matrix.

Taking $rb=b$ into account, we see that $\langle b,p_i\rangle$ is independent of $i$, hence,
$\langle f_0,p_0\rangle=\langle f_0,p_1\rangle=0$. Therefore, $g_1<1$ implies
$\langle f_0,f_0\rangle=(1-g_1)\langle f_0,b\rangle=(1-g_1)^2+2(1-g_1)\langle b,p_0\rangle^2>0$, i.e.,
$f_0\in\Bbb H_\Bbb R^4\cap C_0$. Finally,
$$r\sigma f_0=rf_0-\frac{\langle f_0,p_{n-1}\rangle}{g_2+1}(p_2-p_0)=(1-g_1)b+\langle b,p_0\rangle(p_1+p_2)-$$
$$-\frac{(1-g_1)\langle b,p_0\rangle+\langle b,p_0\rangle(g_1+g_2)}{g_2+1}(p_2-p_0)=(1-g_1)b+\langle
b,p_0\rangle(p_0+p_1)=f_0\ _\blacksquare$$

\medskip

{\bf3.8.~Lemma.} {\sl The isometry\/ $r\sigma|_{C_0}$ of\/ $C_0$ is a rotation by\/ $a$ about\/ $f_0$, where
$$\cos a=\frac{(1-c_1^2)c_mx_1+c_1(1-c_m^2)x_2}{(1-c_1^2)x_1+(1-c_m^2)x_2}.$$}

{\bf Proof.} The isometry $r\sigma$ preserves orientation. By Lemma~3.5, the isometry $r\sigma|_U$ preserves orientation.
Consequently, the isometry $r\sigma|_{C_0}$ also preserves orientation. By Lemma 3.5, it has to be a rotation by some angle
$a$ about $f_0$. Since $\tr(r\sigma|_U)=2g_1$ by Lemma 3.5 and
$$\tr(r\sigma)=\tr r+2g_1+\frac{\langle p_2-p_0,p_1-p_{n-1}\rangle}{g_2+1}=1+2c_1+2c_m+2g_1+\frac{g_1-g_3}{g_2+1}$$
by (3.7), we obtain
$$\cos a=c_1+c_m+\frac{g_1-g_3}{2(g_2+1)}=c_1+c_m+\frac{(c_3-c_1)x_1+(c_{3m}-c_m)x_2}{2(1-c_2)x_1+2(1-c_{2m})x_2}.$$
Taking $c_2=2c_1^2-1$, $c_3=4c_1^3-3c_1$, $c_{2m}=2c_m^2-1$, and $c_{3m}=4c_m^3-3c_m$ into account, we get
$$\cos a=c_1+c_m+\frac{(c_1^3-c_1)x_1+(c_m^3-c_m)x_2}{(1-c_1^2)x_1+(1-c_m^2)x_2}=
\frac{(1-c_1^2)c_mx_1+c_1(1-c_m^2)x_2}{(1-c_1^2)x_1+(1-c_m^2)x_2}\ _\blacksquare$$

{\bf3.9.~Theorem.} {\sl Suppose that the solution of the system
$$\cases(1-c_1)x_1+(1-c_m)x_2=1\\(1-c_1^2)(c_k-c_m)x_1=(1-c_m^2)(c_1-c_k)x_2\endcases$$
satisfies the inequalities\/ $2<(1-c_i)x_1+(1-c_{mi})x_2$ for all\/ $2\leqslant i\leqslant\frac n2$, where\/
$k,m,n\in\Bbb N$, $n$ is even, $1<k<m<\frac n2$, and\/ $c_i:=\cos\frac{2i\pi}n$. Then there exists a disc bundle\/ $M\to S$
over a closed connected orientable surface\/ $S$ admitting a complete real hyperbolic geometry such that\/
$|eM/\chi S|=\frac{4m-4k}{n-4}$.}

\medskip

{\bf Proof.} First, we observe that the solution clearly satisfies the inequalities $0<x_1$, $0<x_2$, and\break
$(1-c_1)x_1+(1-c_m)x_2<2$. The inequality $1<x_1+x_2$ follows from the inequality $2<(1-c_i)x_1+(1-c_{mi})x_2$ with
$i=\frac n2$ because $c_i=-1$ and $c_{mi}=\pm1$. In other words, we get a point $(x_1,x_2)\in R$ in the region $R$.

Let $p(t):=\big(x_1(t),x_2(t)\big)$, $t\in[0,1]$, be a linearly parameterized path in $R$ that joins a point in the segment
$A\subset R$ (see Lemma 3.1) with the point $p(1)=(x_1,x_2)$. The function $a(t)$ is continuous and, by~Lemma 3.8, the
function $\cos a(t)$ has a form $\cos a(t)=\frac{a_1t+a_2}{a_3t+a_4}$ for some $a_1,a_2,a_3,a_4\in\Bbb R$. By Lemma~3.8,
$\cos a(0)=\frac{(1-c_1^2)c_mx_1(0)}{(1-c_1^2)x_1(0)}=c_m$ and
$\cos a(1)=\frac{(1-c_1^2)c_mx_1+c_1(1-c_m^2)x_2}{(1-c_1^2)x_1+(1-c_m^2)x_2}=c_k$ due to the second equation of the system.
Hence, the function $\cos a(t)$ is not constant and is therefore monotonic. Consequently, the~function $a(t)$ is monotonic.

It was understood in Section 2 that the polyhedron $P(t)$ is fibred iff $r\sigma|_{C_0}$ is a periodic isometry whose order
divides $n$, i.e., iff $\cos a(t)=c_j$ for some $j\in\Bbb Z$. Since $a(t)$ is monotonic, $\cos a(0)=c_m$, and
$\cos a(1)=c_k$, we conclude that $eP(1)=m-k$.

By Lemma 3.3, the polyhedron $P(1)$ satisfies the conditions of PPT due to the first equation of the system. It remains to
observe that the Euler characteristic of the corresponding orbifold $S$ equals $\chi S=\frac n4-\frac n2+1=-\frac{n-4}4$
$_\blacksquare$

\medskip

{\bf3.10.~Calculation.} Taking $(k,m,n):=(2,5,24)$ (or $(k,m,n):=(3,6,24)$ or $(k,m,n):=(5,11,44)$), one can check the $11$
(or $11$ or $21$) inequalities of Theorem 3.9 thus arriving at $|eM/\chi S|=\frac35$
$_\blacksquare$

\bigskip

\centerline{\bf References}

\medskip

[AGG] S.~Anan$'$in, C.~H.~Grossi, N.~Gusevskii, {\it Complex hyperbolic structures on disc bundles over surfaces,}
Int.~Math.~Res.~Not.~{\bf2011} (2011), no.~19, 4295--4375, http://arxiv.org/abs/math/0511741

\smallskip

[AGr1] S.~Anan$'$in, C.~H.~Grossi, {\it Yet another Poincar\'e's polyhedron theorem,} Proc.~Edinburgh Math.~Soc., {\bf54}
(2011), 297--308, http://arxiv.org/abs/math/0812.4161

\smallskip

[AGr2] S.~Anan$'$in, C.~H.~Grossi, {\it Coordinate-free classic geometries,} Moscow Math.~J.~{\bf11} (2011), no.~4,
633--655, http://arxiv.org/abs/math/0702714

\smallskip

[AGr3] S.~Anan$'$in, C.~H.~Grossi, {\it Differential geometry of grassmannians and the Pl\"ucker map,} Central European
J.~of Math.~{\bf10} (2012), no.~3, 873--884, http://arxiv.org/abs/0907.4470

\smallskip

[BoS] F.~Bonahon, L.~Siebenmann, {\it The classification of Seifert fibered\/ $3$-orbifolds,} in {\sl Low Dimensional
Topology,} edited by R.~Fenn, LMS Lecture Notes in Science {\bf95}, New York: Cambridge University Press, 1985, 258 pp.

\smallskip

[GLT] M.~Gromov, H.~B.~Lawson Jr., W.~Thurston, {\it Hyperbolic\/ $4$-manifolds and conformally flat\/ $3$-manifolds,}
Inst.~Hautes \'Etudes Sci.~Publ.~Math., no.~68 (1988), 27--45

\smallskip

[Kap] M.~Kapovich, {\it On hyperbolic\/ $4$-manifolds fibered over surfaces,} preprint (1993)\newline
http://www.math.ucdavis.edu/$\widetilde{\phantom{a}}$kapovich/eprints.html

\smallskip

[Kui1] N.~H.~Kuiper, {\it Hyperbolic\/ $4$-manifolds and tessellations,} Inst.~Hautes \'Etudes Sci.~Publ.~Math., no.~68
(1988), 47--76

\smallskip

[Kui2] N.~H.~Kuiper, {\it Fairly symmetric hyperbolic manifolds,} in {\sl Geometry and Topology of Submanifolds,} II
(1990), World Sci.~Publ., 165--204

\smallskip

[Luo] F.~Luo, {\it Constructing conformally flat structures on some Seifert fibred\/ $3$-manifolds,} Math.~Ann. {\bf294}
(1992), no.~3, 449--458

\enddocument